\documentclass {article}

\usepackage {amssymb}
\usepackage {amsfonts}
\usepackage {amsmath}

\newtheorem {thm} {Theorem}
\newtheorem {prop} [thm] {Proposition}
\newtheorem {lemma} [thm] {Lemma}

\newtheorem {defn} {Definition}
\newtheorem {rmk} {Remark}

\begin {document}

\title {An End to End Gluing Construction for Metrics of Constant
Positive Scalar Curvature}
\author {Jesse Ratzkin}

\maketitle

\section {Introduction}


The goal of this paper to to describe a general process by which one
can glue together metrics of constant positive scalar curvature on  
punctured spheres along their ends to obtain new metrics of constant
positive scalar curvature. 

First let $(M_1 = S^n \backslash \{p_0 \dots p_{k_1-1}\}, g_1)$ and
$(M_2 = S^n \backslash \{q_0 \dots q_{k_2-1}\}, g_2)$ be complete
metrics of scalar curvature $n(n-1)$. These metrics are asymptotic to
Delaunay metrics in small (standard spherical) punctured balls about
$p_j$ and $q_j$ respectively. We will refer to these punctured balls
as the ends of $M_1$ and $M_2$. The Delaunay metrics can be written as  
$$u_{\epsilon}^\frac{4}{n-2}(t + T)(dt^2 + d\theta^2)$$ where
$u_\epsilon$ is a periodic function which assumes its minimal value
$\epsilon$ (called the necksize of the metric) at $t=0$. These metrics
are uniquely determined by their singular set on $S^n$, necksize and the
translation parameter $T$. Assume that we choose $g_1$ and $g_2$ such
that the asymptotic necksize of $g_1$ at $p_0$ is equal to the
asymptotic necksize of $g_2$ at $q_0$; we will call this common
necksize $\epsilon$. Then we can truncate $M_1$ and $M_2$ by  removing
small spherical balls around $p_0$ and $q_0$ and patch these two
metrics together at a neck (value of $t$ where $u_\epsilon$ achieves
its minimum) to obtain a new metric $\tilde g$ on $M = S^n \backslash
\{p_1 \dots p_{k_1-1} , q_1\dots q_{k_2-1} \}$. This construction
depends on two parameters: $R$, which we can think of as the size of
the balls we excised in the truncation process, and $\phi \in SO(n)$
which specifies a rotation in the $S^{n-1}$ factor of the second
summand. Notice that the parameter $R$ is discrete. To indicate the
dependence of $\tilde g$ on these parameters we will denote it as
$\tilde g_{R,\phi}$. Much of the analysis is independent of at least
one of these parameters, and in this case we will suppress the
appropriate subscript. We will construct this metric, which we will
call an approximate solution (because its scalar curvature is very
close to $n(n-1)$), in section \ref {approx_soln}.  

The metric $\tilde g_{R,\phi}$ does not have constant scalar
curvature, but the deviation $\psi = n(n-1) - \tilde S_{R,\phi}$ is
globally small. More precisely, without any modification to $g_1$ and
$g_2$, $\psi$ is compactly supported and 
$\|\psi\|_{C^{0,\alpha}(M)} = O(e^{-R})$. After modifying $g_1$ and
$g_2$ by conformal transformations, we can further arrange that
$\|\psi\|_{C^{0,\alpha}(M)} = O(e^{-\gamma_{n+1}(\epsilon) R})$ where
$\gamma_{n+1} (\epsilon)$ is a coefficient we will discuss in section
\ref{delaunay}. We wish to deform $\tilde g_{R,\phi}$ by a 
conformal factor to obtain a metric with scalar curvature
$n(n-1)$. Recall how the scalar curvature transforms under a conformal
change of metric: if $g' = u^{\frac{4}{n-2}}g$ then 
$$S_{g'} = S_g u^{-\frac{4}{n-2}} - \frac{4(n-1)}{n-2}
u^{-\frac{n+2}{n-2}} \Delta_g u ,$$
which we can rewrite as 
\begin {equation} \label {scal_curv_eqn}
\Delta_g u - \frac{n-2}{4(n-1)} S_g u + \frac{n-2}{4(n-1)} S_{g'}  
u^{\frac{n+2}{n-2}} = 0. \end {equation}
If we normalize the scalar curvatures by setting $S_g = n(n-1) - \psi$
and $S_{g'} = n(n-1)$ the above equation becomes
$$\Delta_g u  - \frac{n(n-2)}{4} u + \left ( \frac{n-2}{n-1} \right )
(\psi/4) u+ \frac{n(n-2)}{4} u^{\frac{n+2}{n-2}} = 0.$$
The linearized equation (linearized about $u=1$) is
\begin {equation} \label {jacobi_eqn}
L_g(u) = \Delta_g u + nu + \left ( \frac{n-2}{n-1} \right ) (\psi/4)
u= 0. \end {equation} 
We will call $L_g$ the Jacobi operator associated to $g$ and solutions
$L_g u = 0$ Jacobi fields of $g$. Notice that $L_g = \Delta_g + n$ if
$\psi = 0$, i.e. if we wish to deform a constant scalar curvature
metric into another conformal constant scalar curvature metric. In
the present case $L_{\tilde g_{R,\phi}}$ is a small perturbation of
$\Delta_{\tilde g_{R,\phi}} + n$. 

Our first step is to show that we can solve this linearized equation,
with uniform (in $R$) estimates on the size of the solution
operator. One can also think of this step as finding a positive lower bound
on the spectrum of $L_{\tilde g_{R,\phi}}$ as an operator
between appropriate function spaces. We will address precisely which
function spaces are the proper ones for this problem in section
\ref{gen_lin_anal}. In order to 
prove this we will need to assume that $(M_i, g_i)$ are both
unmarked nondegenerate; i.e. there are no Jacobi fields which decay at
a rate faster than $e^{-t_j}$ on all ends. We will also
need to assume that one can 
adjust the necksize of the end corresponding to $p_0$ in the moduli
space of constant scalar curvature metrics (see the statement of the
theorem for a precise statement of this condition). We need this last
condition to exclude certain Jacobi fields which we could glue
together to yield an exponentially small eigenvalue. This linear
analysis will occupy section \ref{lin_anal}. 

Then in section \ref{deform} we will
explicitly write down geometric deformations of the metric associated
to a parameter $u \in \mathbb{R}^{k(2n+2)}$, where $k = k_1 + k_2 -2$ is
the number of ends of $M$. One can think of these deformations as
adjusting the necksizes and position of the necks of the metric on the
ends of $M$. Finally, in section \ref{contract} we will use the
solution operator we found and these geometric deformations to solve
the nonlinear problem via the Contraction Mapping principle. This will 
yield the existence part of the following theorem. 

\begin {thm} \label {main-thm}
Let $(M_1 = S^n \backslash \{p_0 \dots p_{k_1-1}\}, g_1)$ and
$(M_2 = S^n \backslash \{q_0, \dots q_{k_2-1}\}, g_2)$ be complete
metrics with scalar curvature $n(n-1)$. Assume $g_1$ and $g_2$ are
unmarked nondegenerate and that the asymptotic necksizes associated to $p_0$
and $q_0$ are both $\epsilon$. Assume also that there is a
one-parameter family of scalar 
curvature $n(n-1)$ metrics $g_t$ on $M_1$, $t \in (-\delta, \delta)$, 
where the asymptotic necksize of $g_t$ associated to $p_0$ is $\epsilon
+ t$. Then for $\eta >0$ there is an $R_0$ such that for $R \geq R_0$
one can deform the approximate solution $(M_R, \tilde g_R)$ first by a
geometric parameter $u$ with $|u|< \eta$ and then by a conformal
factor $(1+v)^{\frac{4}{n-2}}$ (with $v$  exponentially decaying) to
obtain a metric with scalar curvature $n(n-1)$. Moreover, this metric
is unmarked nondegenerate. 
\end {thm}

We first remark that the connect sum of two Delaunay metrics
constructed in \cite {MPU1} and all the metrics constructed by Byde in
\cite {B} and by Mazzeo and Pacard in \cite {MP} satisfy all the
hypotheses of this theorem. One particular application of this gluing
construction is to take $(M_1, g_1)$ and $(M_2, g_2)$ to be isometric
and attach $M_1$ to $M_2$ along isometric ends. We will call this
construction doubling along an end. 

One can think of this theorem as the scalar curvature analogue of a
similar end to end gluing construction for surfaces of constant mean
curvature in Euclidean space (see \cite {R} and \cite {MPPR}). In
fact, most of the analysis is the same for the two constructions. This
phenomenon has been widely noted (compare, e.g., \cite {MPU2} and
\cite {KMP}), but not completely explained. 

This theorem is also very much in the spirit of the results of Schoen
in \cite {S}, of Mazzeo, Pollack and Uhlenbeck in \cite {MPU1} and of
Joyce in \cite{J}. In all cases one constructs an approximate solution 
to the gluing problem, solves the linearized equation with uniform
estimates, and solves the nonlinear problem using a fixed point
theorem or an iteration technique. The constructions of Mazzeo and
Pacard in \cite {MP} and of Byde in \cite{B} are similar in spirit,
but use a different method, in that they solve boundary value problems
on appropriate subdomains and then match Cauchy data. 

I would like to thank F. Pacard for suggesting this problem. I
would also like to thank D. Pollack, R. Mazzeo and F. Pacard
for many useful suggestions as I was learning this subject.

\section {Notation} \label {notation}

In this section we establish some notation for the rest of the paper. 

\subsection {Notation for Delaunay Metrics} \label {del_notation}

First we consider the Delaunay metrics. These can be written as 
$$g_{\epsilon,T} = u_\epsilon ^\frac{4}{n-2} (t+T)(dt^2 + d\theta^2)$$
on $\mathbb{R} \times S^{n-1}$. In the case $T=0$ we will suppress it
from the subscript. The function $u_\epsilon$ satisfies the ordinary
differential equation  
$$u'' - \frac{(n-2)^2}{4} u + \frac{n(n-2)}{4} u^{\frac{n+2}{n-2}} =
0.$$ From this ODE one can show that $u_\epsilon$ is a periodic
function uniquely determined by its minimal value $\epsilon$ (once we
normalize $u_\epsilon$ so it achieves its minimum at $t=0$). We will
denote the period of $u_\epsilon$ by $T_\epsilon$. As we will see in
section \ref{delaunay}, solutions of $L_{g_\epsilon}(u) = 0$ which lie
outside a specific two-dimensional space satisfy a bound which (up to
the change of variables $t \mapsto -t$) we can state as 
$$|u(t,\theta)| \left \{ \begin {array} {rcl} 
\leq & ce^t & t<0 \\ 
\geq & c e^{-t} & t>0. \end{array} \right.$$
For more discussion about the solutions to $L_{g_\epsilon}u = 0$, see
section \ref{delaunay}.

\subsection {Notation for Everything Else in this Paper} 
\label {else_notation}  

Recall that we are starting with $(M_1 = S^n \backslash \{p_0, \dots,
p_{k_1-1} \}, g_1)$ and $(M_2 = S^n \backslash \{q_0, \dots ,
q_{k_2-1} \}, g_2)$ two complete metrics with scalar curvature
$n(n-1)$. We will assume that the $p_j$ are mutually disjoint and that
the $q_j$ are mutually disjoint (we will allow, however, $p_j =
q_{j'}$ for some $j$ and $j'$). Let $r_0$ be small enough so that the
discs $B_{r_0}(p_j)$ in the usual round metric are pairwise disjoint,
and also so that the discs $B_{r_0}(q_j)$ in the usual round metric
are pairwise disjoint. Let $M_1^c = S^n \backslash (\cup
B_{r_0}(p_j))$ and $M_2^c = S^n \backslash (\cup B_{r_0}(q_j))$. Next
fix two cutoff functions $\chi_1$ and $\chi_2$ such that 
$$\chi_1 (p) = \left \{ \begin {array} {rl} 0 & p \in M_1^c \\
1 & p \in B_{r_0/2}(p_j)\backslash \{p_j\} \end {array} \right. $$
and 
$$\chi_2(p) = \left \{ \begin {array} {rl} 0 & p \in M_2^c \\ 
1 & p \in B_{r_0/2}(q_j) \backslash \{q_j\}. \end {array} \right.$$

Inside $B_{r_0}(p_j)$ let $r_j(p)$ be the distance in the spherical
metric to $p_j$ and let $t_j = -\log (r_j/ r_0)$. Similarly, in
$B_{r_0}(q_j)$ let $\rho_j(p)$ be the distance is the spherical metric
to $q_j$ and let $\tau_j = -\log (\rho_j/ r_0)$. Then with respect to these
coordinates the asymptotics theorem (see \cite {CGS}  or \cite {KMPS})
states that we can write the metric $g_1$ in $B_{r_0}(p_j)$ as  
$$g_1 = (u_{1, j} + u_{\epsilon_j}(\cdot + T_j))^\frac{4}{n-2} (t_j,
\theta_j) (dt_j^2 + d\theta_j^2)$$ where 
$$\|u_{1,j}\|_{C^{2,\alpha}((\hat t_j -1, \hat t_j + 1) \times
S^{n-1})} = O(e^{- \hat t})$$ 
for $\hat t_j \geq 1$. Similarly, with respect to the
coordinates $(\tau_j, \theta_j)$, one can write 
$$g_2 = (u_{2,j} + u_{\epsilon_j'} (\cdot + T_j'))^\frac{4}{n-2}
(\tau_j, \theta_j)(d\tau_j^2 + d\theta_j^2)$$
where 
$$\|u_{2,j}\|_{C^{2,\alpha} ((\hat \tau_j-1 ,\hat \tau_j+1) \times
S^{n-1})} = O(e^{- \hat \tau_j})$$ for $\hat \tau_j
\geq 1$. We will see later that we can improve these estimates on
$u_{1,0}$ and $u_{2,0}$ using conformal transformations of $S^n$. We
will assume that $\epsilon_0 = \epsilon_0' = \epsilon$ is fixed
throughout the rest of the paper.   

\section {Delaunay Metrics} \label {delaunay}

In this section we will discuss some of the important features of the
$S^{n-1}$ invariant, complete, scalar curvature $n(n-1)$ metrics on
$\mathbb{R} \times S^{n-1}$, which are known as Delaunay metrics. Most
importantly, we discuss the spectral behavior of the operator
$L_{g_\epsilon} = \Delta_{g_\epsilon} + n$ where $g_\epsilon$ is a
Delaunay metric. 

Recall that we can write the Delaunay metrics as 
$$g_\epsilon = u_\epsilon^{\frac{4}{n-2}} (t) (dt^2 + d\theta^2)$$
where $u_\epsilon$ solves the ordinary differential equation 
$$u'' - \frac{(n-2)^2}{4} u + \frac{n(n-2)}{4} u^{\frac{n+2}{n-2}} =
0.$$ 
We remark that solutions to this ODE exist for all time because the
equation has a conserved energy 
$$H = (u_\epsilon ')^2 - \frac{(n-2)^2}{4}u_\epsilon^2 +
\frac{(n-2)^2}{4} u_\epsilon^{\frac{2n}{n-2}},$$ which would become
unbounded if $u_\epsilon$ were to become unbounded. 
These metrics are uniquely determined by their singular set on $S^n$,
the minimum value $\epsilon$ of the conformal factor and a translation
parameter $T$. In 
fact, varying either parameter yields a one-parameter family of
Delaunay metrics. Taking the derivative of this one parameter family,
we obtain two linearly independent solutions to the Jacobi equation:
$$(\Delta_{g_\epsilon} + n) v_\epsilon^{0,\pm} = 0.$$ 
More precisely, we can write 
$$v_\epsilon^{0,+} = \left. \frac{d}{dT}\right|_{T=0} u_\epsilon (\cdot
+ T) \qquad v_\epsilon^{0,-} = \left. \frac{d}{d\eta} \right |_{\eta =
  0} u_{\epsilon + \eta}.$$
From this
construction we see that $v_\epsilon^{0,\pm}$ are independent of
$\theta$, and thus they satisfy the ordinary differential equation 
$$(v_\epsilon^{0,\pm}) '' + 2 \frac{u_{\epsilon}'}{u_\epsilon}
(v_\epsilon^{0,\pm})' + n u_\epsilon^{\frac{4}{n-2}}
v_\epsilon^{0,\pm} = 0.$$ We will normalize 
$v_\epsilon^{0,\pm}$ by choosing the initial conditions 
$$v_\epsilon^{0,+}(0) = 1 \quad v_\epsilon^{0,-}(0) = 0 \quad 
(v_\epsilon^{0,+})'(0) = 0 \quad (v_\epsilon^{0,-})'(0) = 1.$$ 

Indeed, if we try to separate variables for a general solution $v$ of 
$$\Delta_{g_\epsilon} v + n v = 0,$$ we find that $$v(t,\theta) = \sum
v_j(t) \eta_j(\theta)$$ where $\eta_j$ is the $j$th eigenfunction of
$\Delta_{S^{n-1}}$ with eigenvalue $\lambda_j$ (counted with
multiplicity) and $v_j$ satisfies the ordinary differential equation 
\begin {equation} \label {eigen-eqn}
v_j '' + 2\frac{u_{\epsilon}'}{u_\epsilon}v_j'+
(n-\lambda_j)u_\epsilon^{\frac{4}{n-2}} v_j = 0.\end {equation} 
The functions $v_\epsilon^{0,\pm}$ form a basis for the solution space
to this ODE when $j=0$. We will again choose a normalized pair of
solutions to this ODE, $v_\epsilon^{j,+}$ and $v_\epsilon^{j,-}$,
normalized so that 
$$v_\epsilon^{j,+}(0) = 1 \quad v_\epsilon^{j,-}(0) = 0 \quad
(v_\epsilon^{j,+})'(0) = 0 \quad (v_\epsilon^{j,-})'(0) = 1.$$ 

In fact, we can also find the a basis for the solution space when $j =
1, \dots, n$, again by taking explicit geometric deformations of the
metric. To find these deformations, we use stereographic projection to
write the Delaunay metric as  
$$\hat u_\epsilon^{\frac{4}{n-2}} dx^2$$ where $dx^2$ is the standard
Euclidean metric on $\mathbb{R}^n$ and $u>0$ has a singularity at the
origin. Now we can deform this metric by taking translates $a \mapsto
\hat u_\epsilon (\cdot + a)$, so we obtain a Jacobi field by pulling back 
$$\left. \frac{d}{da_j} \right |_{a = 0}\hat u_\epsilon(\cdot + a) =
e^{-t} (-\hat v_\epsilon^{0,+}(t)  + \frac{n-2}{2} \hat u_\epsilon(t)) 
\eta_j(\theta)$$ where $\eta_j$ is the $j$th eigenfunction of
$\Delta_{S^{n-1}}$ (see \cite{KMPS}). Notice in particular that these
functions all decay like $e^{-t}$. To find the other solution to
equation (\ref{eigen-eqn}), we first take the Kelvin transform of $\hat 
u_\epsilon(x) \mapsto |x|^{2-n} \hat u_\epsilon(x/|x|^2)$, translate
as before and then take the Kelvin transform again. One can think of
this deformation as a translation at infinity. Also, one can show that
the Jacobi field associated to this deformation grows like $e^t$
(again, see \cite {KMPS}). 

At this point we introduce the {\it indicial roots} of $L_{g_\epsilon} =
\Delta_{g_\epsilon} + n$, denoted $\gamma_j(\epsilon)$. These are the
exponential growth rates of the solutions to equation
(\ref{eigen-eqn}). From the above computations, one sees that
$\gamma_0(\epsilon) = 0$ and $\gamma_j(\epsilon) = 1$ for $j = 1,
\dots, n$. By the maximum principle, $\gamma_j(\epsilon) > 1$ for $j >
n$. Indeed, $\lambda_j \geq 2n$ for $j \geq n+1$, so the zero order
term in this ODE is $(n-\lambda_j) \hat u_\epsilon^\frac{4}{n-2} \leq
-n$. It is rather remarkable that one can compute $\gamma_j(\epsilon)$
for $j = 0, \dots, n$ and that they are independent of $\epsilon$, but
the other indicial roots are quite hard to compute and probably depend
on $\epsilon$ in some nontrivial way. 

Another way to recover the indicial roots is to conjugate the operator
$\Delta_{g_\epsilon} + n$ by an exponential function $e^{\delta t}$
and the Fourier-Laplace transform. Then one obtains a one (complex)
parameter family of operators on a fixed function space, which varies
analytically with the parameter. By the Analytic Fredholm Theorem,
this family of operators has a meromorphic solution operator, and the
indicial roots turn out to be the imaginary parts of the poles of this
solution operator. In fact, they show that any solution $u$ to 
$$L_{g_\epsilon} u = 0$$ has an asymptotic expansion 
$$u(t,\theta) \sim \sum_{j \geq 0} (a_{j,+}v_\epsilon^{j,+}(t) + a_{j,-}
v_\epsilon^{j,-}(t)) \eta_j(\theta)$$ where $\eta_j$ is the
$j$th eigenfunction of $\Delta_{S^{n-1}}$ (counting multiplicity) and
$v_j$ are the particular solution of equation (\ref{eigen-eqn}) listed
above. In particular,  
$$|v_\epsilon^{j,\pm} (t)| = O(e^{\pm \gamma_j(\epsilon) t}).$$ See
\cite {MPU2} for more about this approach.

A more thorough explanation of the indicial roots in the mean curvature
setting occur in \cite {MPPR}, including an explanation of why they
are called ``indicial roots.''

To sum up this discussion:
\begin {itemize} 
\item $v_\epsilon^{0,+}$ is bounded and periodic and arises from
  translating the neck of the Delaunay metric towards the singularity
\item $v_\epsilon^{0,-}$ is linearly growing and arises from changing
  the necksize of the Delaunay metric
\item $v_\epsilon^{j, \pm}$ grow/decay like $e^{\mp t}$ for $j = 1 ,
  \dots, n$ and both arise from translating the singular set of the
  Delaunay metric
\item in fact, any solution $u$ to $L_{g_\epsilon} u = 0$ which is
  $L^2$ orthogonal to $v_\epsilon^{0,\pm}$ on $S^{n-1}$-cross-sections
  has an   expansion $u(t,\theta) \sim \sum_{j\geq 1} v_j(t)
  \eta_j(\theta)$ where $|v_j(t)| = O(e^{\pm\gamma_j(\epsilon)t})$
  with $\gamma_1(\epsilon) = \gamma_2(\epsilon) = \cdots = 
  \gamma_n(\epsilon) = 1$ and $1 < \gamma_{n+1}(\epsilon) \leq
  \gamma_{n+2} (\epsilon) \leq \cdots \rightarrow \infty$ (we have to
  exclude the $v_\epsilon^{0,\pm}$ terms because one of them grows
  linearly). 
\end {itemize}
One can find rigorous proofs of the above facts in \cite{MPU2},
\cite{MP} and \cite {KMPS}.

\section {The Approximate Solution} \label {approx_soln}

In this section we construct the approximate solution $\tilde
g_{R,\phi}$. 

First we choose some $R= m T_{\epsilon_0}$ for some positive integer
$m$ and $\phi \in SO(n)$ and define $M$ by 
$$M = (M_1 \backslash B_{r_0e^{-(T_0 + R + 1)}} (p_0)) \cup (M_2
\backslash B_{r_0e^{-(T_0' + R + 1)}} (q_0)) / \sim$$ where we
identify $(t_0, \theta_0)$ with $(\tau_0, \phi \theta_0)$ if $t_0 \geq
T_0 + R -1$ and $\tau_0 \geq T_0' + R -1$ and $t_0 + \tau_0 = T_0 +
T_0' + 2R$. The balls $B_r(p_0)$ and $B_r(q_0)$ are balls in the
standard round metric. We will let $C_R$ be the cylinder $\{(t_0,
\theta_0) : T_0 + R - 1 \leq t_0 \leq T_0 + R + 1 \} \sim \{(\tau_0,
\theta_0) : T_0' + R - 1 \leq \tau_0 \leq T_0' + R + 1 \}$. We will
also find it convenient in the following sections to define the
extended cylinder  
$$\hat C_R = (B_{r_0 e^{-T_0}}(p_0) \backslash B_{r_0 e^{-T_0 - R +
1}}(p_0)) \cup C_R \cup (B_{r_0 e^{-T_0'}}(q_0) \backslash B_{r_0
e^{-T_0' - R + 1}}(q_0)),$$ parameterized by $(t, \theta) \in [-R, R]
\times S^{n-1}$. The relationship between $t$ and $t_0$ or $\tau_0$ is
given by $t = t_0 - R - T_0$ for $t < 0$ and $t = -\tau_0 + R + T_0'$
for $t > 0$. This relationship for between $t$ and $t_0$ and $\tau_0$
agrees with the identification of $(t_0, \theta)$ with $(\tau_0,
\phi\theta)$ in $C_R$ listed above. 

Now we will define the metric $\tilde g_{R,\phi}$. First pick a
cutoff function $\chi$ on $M$ such that 
$$\chi(p) = \left \{ \begin {array} {rl} 1 & p \in M_1 \backslash
B_{r_0e^{-(T_0 + R -1)}}(p_0) \\ 0 & p \in M_2 \backslash
B_{r_0e^{-(T_0' + R -1)}}(q_0). \end {array} \right .$$ 
We define the metric $\tilde g_{R,\phi}$ by letting $\tilde g_{R,\phi}
= g_1$ on $M_1 \backslash B_{r_0e^{-(T_0 + R -1)}}(p_0)$,
letting $\tilde g_{R,\phi} = g_2$ on $M_2 \backslash B_{r_0e^{-(T_0' +
R -1)}}(q_0)$ and by letting  
$$\tilde g_{R,\phi}(t_0, \theta_0) = (u_\epsilon (T_0 + t_0) +
\chi(t_0,\theta_0)u_{1,0}(t_0, \theta_0) + (1-\chi(t_0, \theta_0))
u_{2,0}(T_0' +T_0 +2R- t_0, \phi\theta_0))^\frac{4}{n-2} (dt_0^2 +
d\theta_0^2)$$ on $C_R$. The analysis below will often be independent
of at least one of the parameters $R$ and $\phi$; in this case we will
suppress the appropriate subscripts.

We denote the scalar curvature of $\tilde g_{R,\phi}$ by $\tilde
S_{R,\phi}$. Outside of $C_R$, $\tilde g_{R,\phi}$ is either $g_1$ or
$g_2$, and so $\tilde S_{R,\phi} = n(n-1)$ is these regions. A priori,
we also have $$\|u_{1,0}\|_{C^{2,\alpha} (C_R)} = O(e^{-R}),$$ (and a
similar estimate for $u_{2,0}$). However, we can adjust $g_1$ and
$g_2$ by conformal transformations as follows. The term $u_{1,0}$ has
an asymptotic expansion near $p_0$ as   
\begin {equation} \label {asymp_exp} 
u_{1,0} \sim \sum_{j > 1} (a_{j,+} v_\epsilon^{j,+}(t_0) + a_{j,-}
v_\epsilon^{j,-}(t_0)) \eta_j(\theta_0). \end {equation}
The functions $v_\epsilon^{j,\pm}$ correspond  
explicitly to translations of the origin or infinity once under the 
stereographic projection which sends $p_0$ to infinity. So change $g_1$ by
the conformal motion which translates the origin by $(-a_{1,+}, \dots
, -a_{n,+})$ and then by the conformal motion which translates
infinity by $(-a_{1,-}, \dots, a_{n,-})$. This has the effect of
eliminating the first $n$ terms in the expansion (\ref{asymp_exp}),
and so the new metric, which we will still call $g_1$, has an expansion
of the form $(u_\epsilon + u_{1,0})^\frac{4}{n-2}(dt_0^2 +
d\theta_0^2)$ where now
$$u_{1,0} \sim \sum_{j \geq n+1} (a_{j,+} v_\epsilon^{j,+}(t_0) +
a_{j,-} v_\epsilon^{j,-}(t_0)) \eta_j(\theta_0),$$ and so 
$$\|u_{1,0}\|_{C^{2,\alpha}(C_R)} = O(e^{-\gamma_{n+1} (\epsilon)
  R}).$$ We can perform as similar adjustment to $g_2$ so that
$\|u_{2,0}\|_{C^{2,\alpha}(C_R)} = O(e^{-\gamma_{n+1} (\epsilon)
  R})$. Notice we can only do this adjustment for one end of each of
  the $M_i$. The geometric effect of this adjustment is to translate
  the $p_j$ and $q_j$ (for $j \geq 0$) around so as to make the metrics
  $g_1$ and $g_2$ near $p_0$ and $q_0$ (respectively) closer to being
  Delaunay metrics. The above estimates imply the following lemma.   

\begin {lemma}
When $(M, \tilde g_{R,\phi})$ is defined as above, $\tilde S_{R,\phi}
= n(n-1) - \psi$ where $\psi$ is compactly supported and
$\|\psi\|_{C^{0,\alpha}(M)} = O(e^{-\gamma_{n+1}(\epsilon)R})$.
\end {lemma} 

{\bf Proof}: We only need to estimate $\tilde S_{R,\phi}$ in $C_R$. To
this end, we first rewrite the conformal factor on $C_R$ as
$u_\epsilon(T_0 + t_0) + \chi
u_{1,0}(t_0,\theta_0) + (1-\chi)u_{2,0}(T_0'+T_0+2R- t_0, \phi\theta_0) =
u_\epsilon (T_0+ t_0)(1 + v(t_0,\theta_0))$. If we plug
$1 + v$ into equation (\ref{scal_curv_eqn}), we find 
\begin {eqnarray*}
\left ( \frac{n-2}{n-1} \right ) \frac{\psi}{4}(1+v)^\frac{n+2}{n-2}& =&
\Delta_{g_\epsilon}(1+v) - 
\frac{n(n-2)}{4}(1+v) + \frac{n(n-2)}{4} (1+v)^\frac{n+2}{n-2} \\ 
& = & \Delta_{g_\epsilon} (v) + nv + O(\|v\|_{C^{2,\alpha}}^2).
\end {eqnarray*}
The lemma now follows from the above bounds on $u_{1,0}$ and
$u_{2,0}$, which imply similar bounds on $v$. \hfill $\blacksquare$

One way to rephrase the result of this lemma is to say that one can
write the metric $\tilde g_{R,\phi}$ restricted to $\hat C_R$ (in the
$(t,\theta)$ coordinates) as 
$$\tilde g_{R,\phi} = (u_\epsilon(t) + v(t, \theta))^\frac{4}{n-2}
(dt^2 + d\theta^2)$$ where 
$$|v(t,\theta)| = O\left ( \frac{\cosh^{\gamma_{n+1}(\epsilon)}
  t}{\cosh^{\gamma_{n+1}(\epsilon)} R} \right ) .$$  
Because of the above bounds on $\tilde S_{R,\phi}$ we will call the
metric $\tilde g_{R,\phi}$ an approximate solution to our problem. 

\section {Linear Analysis} \label{lin_anal}

In this section we will develop the necessary linear analysis to find
a uniformly bounded solution operator for the Jacobi operator
$L_{\tilde g_{R,\phi}}$. We start the section by recalling some of the
linear analysis for constant positive scalar curvature metrics in
\cite {MPU2} and then we construct a solution operator for $L_{\tilde
  g_{R,\phi}}$.  

\subsection {Linear Analysis for General Constant Scalar Curvature
Metrics on Punctured Spheres}
\label {gen_lin_anal}

The growth properties for solutions of $L_{g_\epsilon}(u) = 0$
outlined above motivate the use of the following function spaces. 
\begin {defn}
On $(M_i, g_i)$ we define $C^{l, \alpha}_\delta (M_i)$ to
be the space of functions such that the norm 
$$\|u\|_{C^{l,\alpha}_\delta(M_i)} = \|u\|_{C^{l,\alpha}(M_i^c)} +
\max_{0 \leq j \leq k_1-1} \sup_{\hat t_j \geq e^{-r_0} + 1}
\|e^{-\delta t_j} u\|_{C^{l,\alpha}((\hat t_j -1, \hat t_j + 1) \times
S^{n-1})}$$ is finite. There is a similar definition for $(M_2,
g_2)$. For the approximate solution $(M, \tilde g_{R,\phi})$  we will
need to adjust this definition as follows. Recall that we can write 
$$M = M_1^c \cup (\cup_1^{k_1-1} B_{r_0}(p_j) \backslash\{p_j\}) \cup
M_2^c \cup (\cup_1^{k_2-1} B_{r_0}(q_j) \backslash\{q_j\}) \cup \hat
C_R.$$ Then we define $C^{l,\alpha}_\delta(M)$ to be the space of
functions such that the norm 
\begin {eqnarray*}
\|u\|_{C^{l,\alpha}_\delta(M)} & = & \|u\|_{C^{l,\alpha}(M_1^c)} +
\|u\|_{C^{l,\alpha} (M_2^c)} + \max_{1 \leq j \leq k_1-1} \sup_{\hat
t_j \geq e^{-r_0} + 1} \|e^{-\delta t_j} u\|_{C^{l,\alpha}((\hat t_j
-1, \hat t_j + 1) \times S^{n-1})} \\ 
&& + \max_{1 \leq j \leq k_2-1} \sup_{\hat \tau_j \geq
e^{-r_0} + 1} \|e^{-\delta \tau_j} u\|_{C^{l,\alpha}((\hat \tau_j -1,
\hat \tau_j + 1) \times S^{n-1})} 
 + \sup_{|\hat t| \leq R-1} \|\frac{\cosh^\delta
R}{\cosh^\delta t} u\|_{C^{l,\alpha} ([\hat t -1, \hat t + 1] \times
S^{n-1})} \end {eqnarray*} is finite. We 
also say $(M_i, g_i)$ is unmarked nondegenerate if $L_{g_i} u = 0$ does
not admit solutions $u \in C^{2,\alpha}_{\delta}(M_i)$ for any $\delta
< -1$. 
\end {defn}

Functions in $C^{l,\alpha}_\delta (M_i)$ can grow at most like
$e^{\delta t_j}$ on the end $E_j$. We remark that for $\delta \leq 1$
the only solutions of $\Delta_{g_\epsilon}v + nv = 0$ with $v \in
C^{2,\alpha}_\delta (\mathbb{R} \times S^{n+1})$ are linear
combinations of $v_\epsilon^{j,+}$ and $v_\epsilon^{j,-}$ for $j = 0,
\dots, n$, each of which is either bounded and periodic or unbounded
on at least one end. So the Delaunay metrics are unmarked
nondegenerate. We also remark that the function space
$C^{l,\alpha}_\delta (M)$ is the same space of functions as if we had
not weighted the middle cylinder, but it has a  different norm. This
difference in norms will become important later when we want uniform
bounds on a solution operator. We further remark that the notion of
unmarked nondegenerate is weaker than the notion of marked
nondegenerate, which requires that $L_{g_i} u = 0$ does not have any
solutions where  $u \in C^{2,\alpha}_\delta (M_i)$ with $\delta < 0$.  

In order to find a function space on which $L_{g_i}$ has suitable
mapping properties we will need the following definition.
\begin {defn} 
The deficiency space $W_{g_i}$ of $(M_i, g_i)$ is
the span of all the functions $\chi_i v_{\epsilon_j}^{i,\pm}$, where $1 \leq
j \leq k_i$ (so the sum runs over all the ends of $M_i$) and $0 \leq i
\leq n$. 
\end {defn}
Notice that $W_{g_i}$ is a vector space of dimension $k_i (2n +2)$
and it has a basis $\{ \chi_i v_{\epsilon_j}^{i,\pm} \}$ which only
depends on the metric $g_i$ and the choice of cutoff function
$\chi_i$. We will use this basis to give $W_{g_i}$ the Euclidean
norm. We also remark that this is the proper deficiency space to use
to parameterize the unmarked moduli space, meaning that one fixes the
cardinality but not the position of the singular set. For the marked
moduli space (fixing both the cardinality and the position of the
singular set) one should work with a smaller deficiency space which
only incorporates the Jacobi fields arising from translating the necks
of the Delaunay metrics along their axes and changing their necksizes. 

\begin {rmk} \label {translations}
Let $E_j$ be an end of $M_1$ corresponding to the puncture point
$p_j$. It turns out the we can use particular conformal Killing fields
on the sphere to show that for any end $E_j$ there is always a Jacobi
field of $g_i$ which is asymptotic to $v_{\epsilon_j}^{0,+}$ along
$E_j$. To see this, consider stereographic projection sending $p_j$ to
$\infty$ composed with a dilation about the origin. This provides a 
one-parameter family of scalar curvature $n(n-1)$ metrics on $M_1$
which translate the Delaunay neck on $E_j$. Taking the infinitesimal
generator of this family we obtain a Jacobi field asymptotic to
$v_\epsilon^{0,+}$. Similar Jacobi fields exist on $M_2$. These Jacobi
fields are also in $C^{2,\alpha}_1(M_i)$. This seems to
be a special property of spheres, as one cannot in general find such
conformal Killing fields on arbitrary compact manifolds with positive
scalar curvature. In the mean curvature case the corresponding Jacobi
fields arise from global translations of the surface.
\end {rmk}

A similar Linear Decomposition result to the one stated below appears
as Lemma 4.18 of \cite {MPU2}, as stated for weighted Sobolev spaces
and exactly constant scalar curvature metrics. The result below is
essentially the next term in the asymptotic expansion; see Proposition
4.15 of \cite {MPU2}. The proof for weighted H\"older spaces is nearly
identical and really only requires that the ends are asymptotically
Delaunay.  
\begin {prop} \label {lin_decomp}
{\rm (Mazzeo, Pollack, Uhlenbeck, 1996)} Let $\delta \in (1, \inf
\gamma_{n+1}(\epsilon_j))$. If $u \in C^{2,\alpha}_\delta (M_i)$, $f \in
C^{0,\alpha}_{-\delta}(M_i)$ and $L_{g_i} u = f$ then $u \in W_{g_i}
\oplus C^{2,\alpha}_{-\delta}(M_i)$.
\end {prop}

Suppose $g_i$ is a unmarked nondegenerate metric on $M_i$. Then for
$\delta \in (1, \inf (\gamma_{n+1}(\epsilon_j))$ 
$$L_{g_i}: C^{2,\alpha}_{-\delta}(M_i) \rightarrow
C^{0,\alpha}_{-\delta}(M_i)$$ is injective, which in turn implies 
$$L_{g_i} : C^{2,\alpha}_\delta (M_i) \rightarrow
C^{0,\alpha}_\delta(M_i)$$ is surjective. If we combine this with the
Linear Decomposition result in proposition \ref{lin_decomp} then we see
that  
$$L_{g_i} : W_{g_i} \oplus C^{2,\alpha}_{-\delta}(M_i) \rightarrow
C^{0,\alpha}_{-\delta} (M_i)$$ is surjective. We will call the kernel of
this map $B_{g_i}$, the {\em bounded null space} of $L_{g_i}$. Mazzeo,
Pollack and Uhlenbeck (\cite{MPU2}) show that if $M_i$ has $k_i$ ends
and $g_i$ is unmarked nondegenerate then $B_{g_i}$ is $k_i
(n+1)$-dimensional (in general $B_{g_i}$ could contain a space of
exponentially decaying functions of some unknown dimension). From
this reasoning one can see (using the Implicit Function Theorem) that
near an unmarked nondegenerate point the moduli space of such metrics
has the structure of a real analytic manifold of dimension $k_i(n+1)$. 

\subsection {Solvability of the Linear Problem} \label{soln_op}

To construct the deficiency space $W_{\tilde g_{R,\phi}}$ we take
cutoffs of the Jacobi fields $v_\epsilon^{i,\pm}$from the model
Delaunay metrics arising from $p_1 \dots p_{k_1-1}$ and $q_1, \dots
q_{k_2-1}$. Notice we do not include $p_0$ and $q_0$. Again, we will
use the basis formed by $\{\chi_1 v_{\epsilon_j}^{i,\pm}, \chi_2
v_{\epsilon_j'}^{i,\pm}\}$, which induces the Euclidean norm on
$W_{\tilde g_{R,\phi}}$.  

Recall that the Jacobi operator $L_{\tilde g+{R,\phi}}$ is  given by 
$$L_{\tilde g_{R,\phi}} = \Delta_{\tilde g_{R,\phi}} + n + \left (
\frac{n-2}{n-1} \right ) \frac{\psi}{4},$$ which is a perturbation of
$\Delta_{\tilde g_{R,\phi}} + n$ where the perturbation is compactly
supported and globally of size $O(e^{-\gamma_{n+1}(\epsilon)R})$. 

\begin {prop} \label {nondegen}
Suppose both $g_i$ are unmarked nondegenerate and there exists a one-parameter
family of scalar curvature $n(n-1)$ metrics $g_t$ on $M_1$ such that
the asymptotic necksize of the end at $p_0$ with respect to to $g_t$
is $\epsilon + t$. Then for $\delta \in (1,
\inf\{\gamma_{n+1}(\epsilon_j), \gamma_{n+1}(\epsilon'_j)\})$ there exists an
$R_0>0$ such that for $R \geq 0$ one can find an operator  
$$G_{R,\phi} : C^{0,\alpha}_{-\delta} (M)\rightarrow 
W_{\tilde g_{R,\phi}} \oplus C^{2,\alpha}_{-\delta}(M)$$ such that $u
= G_{R,\phi}(f)$ solves the equation $L_{\tilde g_{R,\phi}} (u) = f$
and $\|u\|_{W_{\tilde g_{R,\phi}} \oplus C^{2,\alpha}_{-\delta}(M)}
\leq c\|f\|_{C^{0,\alpha}_{-\delta}(M)}$ where $c$ is independent of
$R$ and $\phi$. 
\end {prop} 

The idea behind this proof was communicated to me by F. Pacard. 

{\bf Proof}: We wish to solve the equation $$L_{\tilde g_{R,\phi}} (u)
= f.$$ To this end, first let $u_1 + v_1 \in W_{g_1} \oplus
C^{2,\alpha}_{-\delta} (M_1)$ solve $$L_{g_1} (u_1 + v_1) = \chi f.$$
Such a solution exists because $g_1$ is unmarked nondegenerate. Moreover, we
have the estimate 
\begin {equation} \label {decay_est1}
\|u_1\|_{W_{g_1}} + \|v_1\|_{C^{2,\alpha}_{-\delta}(M_1)} \leq c_1
\|\chi f\|_{C^{0,\alpha}_{-\delta}(M_1)}. \end {equation} 
In $B_{r_0}(p_0)$ (the standard spherical ball), 
$u_1(t_0,\theta_0) \sim \sum_{i, \pm} \alpha_{i,\pm}
v_\epsilon^{i,\pm}(t_0)$. Now  
choose $\Phi_1 \in B_{g_1}$ such that $|\Phi_1(t_0,\theta_0) + \sum
\alpha_{i,\pm} v_\epsilon^{i,\pm}(t_0)| =
O(e^{-\gamma_{n+1}(\epsilon) t_0})$. We can choose such a $\Phi_1$ because
of the existence of the Jacobi fields in remark \ref{translations} and 
because of the assumption that there is a one-parameter family of
metrics $g_t$ on $M_1$ such that $g_0 = g_1$ and the asymptotic
necksize of the end at $p_0$ with respect to $g_t$ is $\epsilon +
t$. Thus $L_{g_1}(u_1 + v_1 + \Phi_1) = \chi f$ and 
\begin {equation} \label {decay_est2}
|u_1(t_0,\theta_0) + v_1(t_0,\theta_0) + \Phi_1(t_0,\theta_0)| \leq
2c_1 \|\chi f\|_{C^{0,\alpha}_{-\delta} (M_1)} e^{-\delta t_0} 
\end {equation} for $(e^{-t_0}, \theta) \in B_{r_0}(p_0)$. We also
have the estimate 
\begin {equation} \label {decay_est3}
\|\Phi_1\|_{W_{g_1} \oplus C^{2,\alpha}_{-\delta}(M_1)} \leq c_1
\|\chi f\|_{C^{0,\alpha}_{-\delta} (M_1)}. \end {equation}

Similarly we let $u_2 + v_2 \in W_{g_2} \oplus C^{2,\alpha}_{-\delta}
(M_2)$ solve $$L_{g_2}(u_2 + v_2) = (1-\chi) f$$ with the estimate 
\begin {equation} \label {decay_est4}
\|u_2\|_{W_{g_2}} + \|v_2\|_{C^{2,\alpha}_{-\delta} (M_2)} \leq c_2
\|f\|_{C^{0,\alpha}_{-\delta} (M_2)}.\end {equation}
This time we cannot cancel the nondecaying part of $u_2 + v_2$ on
$E_2$. Instead, let $\beta_{i,\pm}$ be such that
$|u_2(\tau_0,\theta_0) - \sum \beta_{i,\pm} v_\epsilon^{i,\pm}
(\tau_0)| = O(e^{-\gamma_{n+1}(\epsilon) \tau_0})$ and let 
$\Phi_2 \in B_{g_1}$ be such that $|\Phi_2(t_0,\theta_0) - \sum
\beta_{i,\pm} v_\epsilon^{i,\pm}(t_0)| =
O(e^{-\gamma_{n+1}(\epsilon)t_0})$ (recall that the relationship between
$t_0$ and $\tau_0$ in $\hat C_R$ is given by $t_0 + \tau_0 = T_0 +
T_0' + 2R$). This time the salient estimates are 
\begin {equation} \label {decay_est5}
\|\Phi_2\|_{W_{g_1} \oplus C^{2,\alpha}_{-\delta}(M_1)} \leq c_2
\|(1-\chi) f\|_{C^{0,\alpha}_{-\delta}(M_2)} \end {equation}
and 
\begin {equation} \label {decay_est6}
\|L_{\tilde g_{R,\phi}}(\Phi_2 + u_2)\|_{C^{0, \alpha} (C_R)} \leq c_3
\|f\|_{C^{0, \alpha}_{-\delta}(M)} e^{-\delta R}. \end {equation}
This last estimate is a straightforward calculation using the facts
that $\Phi_2$ is a Jacobi field of $g_1$, $u_2 \in W_{g_2}$, and both
$g_1$ and $g_2$ (and hence $\tilde g_{R,\phi}$) are
$C^{2,\alpha}$-close to being Delaunay metrics in $C_R$.

Now choose cutoff functions $\eta_1$ and $\eta_2$ such that 
$$\eta_1(p) = \left \{ \begin {array}{rl} 1 & p \in (M_1 \backslash
B_{r_0 e^{-(T_0 + R -1)}}(p_0)) \cup C_R \cup (B_{r_0e^{-(T_0' +2)}}(q_o)
\backslash B_{r_0 e^{-(T_0' + R -1)}}(q_0)) \\ 
0 & p \in M_2 \backslash B_{r_0 e^{-(T_0' + 1)}}(q_0) \end {array}
\right.$$  
and 
$$\eta_2(p) = \left \{ \begin {array}{rl} 1 & p \in (M_2 \backslash
B_{r_0 e^{-(T_0' + R -1)}}(q_0)) \cup C_R \cup (B_{r_0e^{-(T_0 + 2)}}(p_0)
\backslash B_{r_0 e^{-(T_0 + R - 1)}}(p_0)) \\ 0 & p \in M_1 \backslash
B_{r_0 e^{-(T_0 + 1)}}(p_0), \end {array} \right.$$  
and define $$u+v = \hat G_{R,\phi} (f) = \eta_1(u_1 + v_1 +\Phi_1) +
\eta_2 v_2 + (1-\chi) u_2 +\chi \Phi_2.$$ We will complete the proof
of this proposition by showing 
\begin {itemize}
\item $\|u+v\|_{W_{\tilde g_{R,\phi}} \oplus C^{2,\alpha}_{-\delta}
(M)} \leq c \|f\|_{C^{0,\alpha}_{-\delta}(M)}$ and 
\item $\| L_{\tilde g_{R,\phi}}(u+v) - f\|_{C^{0,\alpha}_{-\delta}(M)}
\leq \tilde c \|f\|_{C^{0,\alpha}_{-\delta}(M)} e^{-\tilde \gamma R}$
for some $\tilde \gamma > 0$
\end {itemize}
where $c$ and $\tilde c$ are independent of $R$. The above estimates
show $\hat G_{R,\phi}$ is uniformly bounded and $L_{\tilde g_{R,\phi}}
\circ \hat G_{R,\phi} = \rm{Id} + O(e^{-\tilde \gamma R})$ and the
statement of the proposition follows by a simple perturbation
argument. The first estimate follows immediately from the estimates
(\ref{decay_est1}), (\ref{decay_est3}), (\ref{decay_est4}) and
(\ref{decay_est5}). For the second estimate, notice $L_{\tilde
g_{R,\phi}}(u+v) - f \neq 0$ only where $\nabla \chi$,  $\nabla\eta_1$ or
$\nabla\eta_2$ are nonzero. The region where $\nabla \eta_1 \neq 0$
corresponds to $T_0' + 1 \leq \tau_0 \leq T_0' + 2$, or $R-2 \leq t
\leq R-1$. In this region 
\begin {eqnarray*}
|L_{\tilde g_{R,\phi}}(u+v)(t,\theta) - f(t,\theta) | & = & |L_{\tilde
g_{R,\phi}}(\eta_1(u_1 + v_1 + \Phi_1)(t,\theta)| \\
& \leq & \|L_{\tilde g_{R,\phi}}\| \cdot |u_1(t, \theta) + v_1(t, \theta) +
\Phi_1(t,\theta)|\\ 
& = & O(e^{-2\delta R})\|f\|_{C^{0, \alpha}_{-\delta}(M)}.
\end {eqnarray*}
One can similarly estimate $L_{\tilde g_{R,\phi}}(u+v) -f$ in the
regions $\nabla \chi \neq 0$ and $\nabla \eta_2 \neq 0$.
\hfill $\blacksquare$ 

Whenever we wish to solve the equation $L_{\tilde g_{R,\phi}} (u) = f$
where $f$ decays at some exponential rate we will always use the
solution operator we constructed in the above proposition. In
general, one can find many solution operators for $L_{\tilde
g_{R,\phi}}$ but most of them will not be uniformly bounded in $R$, as
$G_{R,\phi}$ is. 

\section {Nonlinear Analysis}

There are two parts to the nonlinear part of this problem. First we
construct explicit deformations of the metric $\tilde g_{R,\phi}$,
parameterized by a small ball about the origin in $W_{\tilde
g_{R,\phi}}$. These deformed metrics will not be conformal
to $\tilde g_{R,\phi}$, but their conformal class will always lie
close to that of $\tilde g_{R,\phi}$ in the Gromov-Hausdorff
topology (see below). One can think of this step as an exponential map
from a subspace of the tangent space of all metrics to the space of
metrics itself. Finally, we use the solution operator $G_{R,\phi}$ to
to build a contraction from a small ball in $W_{\tilde g_{R,\phi}}
\oplus C^{2,\alpha}_{-\delta} (M)$ to itself. The fixed point of this 
contraction will be the solution to our nonlinear problem. 

\subsection {The Geometric Deformations} \label {deform}

In this section we define geometric deformations of the approximate
solution $\tilde g_{R,\phi}$ corresponding to elements in $W_{\tilde
g_{R,\phi}}$. 

First recall that we can write $w \in W_{\tilde g_{R,\phi}}$ as 
$$w = \chi_1 \sum_{j=1}^{k_1-1} (\alpha_j^{i,+} v_{\epsilon_j}^{i,+} +
\alpha_j^{i,-} v_{\epsilon_j}^{i,-}) \eta_i 
+ \chi_2 \sum _{j=1}^{k_2-1} (\beta_j^{i,+} v_{\epsilon_j'}^{i,+} +
\beta_j^{i,-} v_{\epsilon_j'}^{i,-}) \eta_i$$
where $v_{\epsilon_j}^{i,\pm}$ and 
$v_{\epsilon_j'}^{i,\pm}$ are the Jacobi fields for $g_{\epsilon_j}$ and
$g_{\epsilon_j'}$ described in section \ref{delaunay}. For each 
end of $M_1$, let $\hat g_j$ be the deformed Delaunay
metric obtained as follows: First replace $u_{\epsilon_j}(\cdot +
T_j)$ with $u_{\epsilon_j + \alpha_j^{0,-}}(\cdot + T_j +
\alpha_j^{0,+})$. Then transfer to $\mathbb{R}^n$ via stereographic
projection and translate the origin by $(\alpha_j^{1,+}, \dots,
\alpha_j^{n,+})$ and translate infinity by $(\alpha_j^{1,-}, \dots,
\alpha_j^{n,-})$. Then pull the result back by the inverse of stereographic
projection. The result is a new Delaunay metric $(\hat
u_j)^\frac{4}{n-2}(dt_j^2 + d\theta^2)$. Similarly, we obtain a new
Delaunay metric $(\hat u_{j'})^\frac{4}{n-2}(d\tau_j^2 + d\theta^2)$
on each end of $M_2$. 

Given $w \in W_{\tilde g_{R,\phi}}$ in the above form, we will define
a new metric $\tilde g_{R,\phi}(w)$ as follows. First let $M^c$ be $M
\backslash ((\cup_1^{k_1-1} B_{r_0}(p_j)) \cup (\cup_1^{k_2-1}
B_{r_0}(q_j)))$ where the balls are taken with respect to the standard
spherical metric. Let $\hat \chi$ be a cutoff function such that 
$$\hat \chi(p) = \left \{ \begin {array}{rl} 1 & p \in (\cup_1^{k_1-1}
B_{r_0/2} (p_j) \backslash \{ p_j\}) \cup (\cup_1^{k_2-1}
B_{r_0/2} (q_j) \backslash \{ q_j\}) \\ 0 & p \in M^c. 
\end {array} \right. $$ We will let $\tilde g_{R,\phi}(w) = \tilde
g_{R,\phi}$ on $M^c$. In $B_{r_0}(p_j) \backslash \{p_j\}$ we define
$\tilde g_{R,\phi}$ in the local cylindrical coordinates (see section
\ref{notation}) by 
$$\tilde g_{R,\phi}(w) = (u_{1,j} + (1-\chi)u_{\epsilon_j}(\cdot +
T_j) + \chi \hat u_j(\cdot + T_j))^\frac{4}{n-2}
(dt_j^2 + d\theta_j^2).$$ We define $\tilde g_{R,\phi}(w)$ in
$B_{r_0}(q_j) \backslash \{q_j\}$ similarly. 

In defining the deformed metric $\tilde g_{R,\phi}(w)$ we induce new
perturbations in the scalar curvature, compactly supported in
$B_{r_0}(p_j) \backslash \{p_j\}$ and $B_{r_0} (q_j) \backslash
\{q_j\}$. But this perturbation is small, as one can see by taking a 
Taylor expansion of the scalar curvature as given in equation
(\ref{scal_curv_eqn}). The proof of the following lemma is a
straightforward computation and left to the reader. 
\begin {lemma} \label {small-error}
For $\delta \in (1, \inf \{\gamma_{n+1}(\epsilon_j),
\gamma_{n+1}(\epsilon_j') \})$, the scalar curvature of $\tilde
g_{R,\phi}(w)$ is given by $\tilde S_{R,\phi}(w) = n(n-1) -\psi - \hat
\psi (w)$ where $\|\hat \psi (w)\|_{C^{0,\alpha}_{-\delta} (M)} =
O(|w|)$, but not $o(|w|)$. Moreover, $\hat \psi(w)$ is compactly
supported in $(\cup_{j = 1, \dots, k_1-1} (B_{r_0}(p_j) \backslash
B_{r_0 / 2}(p_j))) \cup (\cup_{j = 1, \dots, k_2-1} (B_{r_0} (q_j)
\backslash B_{r_0/2} (q_j)))$.
\end {lemma} 
Thus in particular the perturbation in the scalar curvature decays
exponentially, so the nonlinear operator is well-behaved. 

Notice that the metric $\tilde g_{R,\phi}(w)$ will not be
conformal to $\tilde g_{R,\phi}$, but for small $w$ the metrics will
be close on large compact sets. More precisely, given a compact set
$\Omega \subset M$ and $\nu > 0$ there is an $\eta >0$ such that
$\|\tilde g_{R,\phi}(w) - \tilde g_{R,\phi}\|_{C^{0,\alpha} (\Omega)}
\leq \nu$ for $|w| \leq \eta$. In other words, $\tilde 
g_{R,\phi}(w)$  and $\tilde g_{R,\phi}$ are close in the
Gromov-Hausdorff topology. This is a more precise way of saying that
the conformal class of $\tilde g_{R,\phi}(w)$ (and hence that of the
constant scalar curvature metric we will construct) is close to the
conformal class of $\tilde g_{R,\phi}$. 

\subsection {Solving the Gluing Problem with a Contraction} \label {contract}

We wish to solve the nonlinear equation \ref{scal_curv_eqn}, which we
will restate here with $S_g = n(n-1) - \psi$ and $S_{g'} = n(n-1)$:
$$0 = \Delta_{\tilde g_{R,\phi}} u - \frac{n(n-2)}{4} u + \left (
\frac{n-2}{n-2} \right ) \frac{\psi}{4} u + \frac{n(n-2)}{4}
u^\frac{n+2}{n-2} = L_{\tilde g_{R,\phi}} u + Q_{R,\phi}(u).$$
Above we have written a Taylor series expansion for the scalar
curvature operator, where $$L_{\tilde g_{R,\phi}} = \Delta_{\tilde
g_{R,\phi}} + n + \left (\frac{n-2}{n-1} \right ) \frac{\psi}{4}$$ and
$$Q_{R,\phi} (u) = \frac{n(n-2)}{4} u^\frac{n+2}{n-2} - 
\frac{n(n+2)}{4} u$$ 
incorporates all the second and higher order terms of the Taylor
series. In particular, 
$$Q_{R,\phi}(0) = 0 \qquad \nabla Q_{R,\phi}(0) = 0,$$
and so $$\|Q_{R,\phi}(u)\|_{C^{0,\alpha}_{-\delta}(M)} \leq C_Q
\|u\|^2_{W_{\tilde g_{R,\phi}} \oplus C^{2,\alpha}_{-\delta}(M)}$$ for
some $C_Q$. 

The above analysis deserves some comment before we continue. Firstly,
we decompose $u$ as $ w + v$ where $w \in W_{\tilde g_{R,\phi}}$ and
$v \in C^{2,\alpha}_{-\delta}(M)$. We make sense of the scalar
curvature operator applied to $u$ as follows: first deform $\tilde
g_{R,\phi}$ to $\tilde g_{R,\phi}(w)$ as in section \ref{deform} and
then let the new metric $g$ be given by $g = (1 + v)^\frac{4}{n-2}
\tilde g_{R,\phi}(w)$. The scalar curvature of the new metric $g$ is
the scalar curvature operator applied to $u$. Finding $u = v+w \in
W_{\tilde g_{R,\phi}} \oplus C^{2,\alpha}_{-\delta}(M)$ such the the
scalar curvature of $g = (1 + v)^\frac{4}{n-2} \tilde g_{R,\phi}(w)$
is $n(n-1)$ is equivalent to solving the equation $L_{\tilde
g_{R,\phi}}(u) = -Q_{R,\phi}(u)$. Moreover, both of these operators
are well defined acting on $W_{\tilde g_{R,\phi}} \oplus
C^{2,\alpha}_{-\delta}(M)$ and map into
$C^{0,\alpha}_{-\delta}(M)$. To see that $Q_{R,\phi} (u)$ decays
exponentially, recall that that the new metric $g$ restricted to the
ends is still conformal to a Delaunay metric with a conformal factor
that is exponentially close to $1$. Notice that $Q_{R,\phi}(u)$
contains terms from the conformal factor $v$ and from the perturbation
term $\hat \psi (w)$. 

The existence part of theorem \ref{main-thm} follows immediately from
the following proposition. 

\begin {prop}
The map $$K_{R,\phi} : W_{\tilde g_{R,\phi}} \oplus
C^{2,\alpha}_{-\delta}(M) \rightarrow W_{\tilde g_{R,\phi}} \oplus
C^{2,\alpha}_{-\delta} (M)$$ given by 
$$K_{R,\phi}(u) = -G_{R,\phi}(Q_{R.\phi}(u))$$
is a contraction on sufficiently a small ball centered at the origin,
and thus it has a unique fixed point. 
\end {prop}

{\bf Proof}: First we estimate 
\begin {eqnarray*}
\|K_{R,\phi}(u)\|_{W_{\tilde g_{R,\phi}} \oplus
C^{2,\alpha}_{-\delta}(M)} & \leq &
\|G_{R,\phi}\|\cdot\|Q_{R,\phi}(u)\|_{C^{0,\alpha}_{-\delta}(M)} \\
& \leq & C_Q \|G_{R,\phi}\| \cdot\|u\|^2_{W_{\tilde g_{R,\phi}} \oplus
C^{2,\alpha}_{-\delta}(M)}, \end {eqnarray*}
which shows that $K_{R,\phi}$ maps a small ball to itself. Also, 
\begin {eqnarray*}
\|K_{R,\phi}(u_1) - K_{R,\phi}(u_2)\|_{W_{\tilde g_{R,\phi}} \oplus
C^{2,\alpha}_{-\delta}(M)} & \leq & C_Q \|G_{R,\phi}\| \cdot\|u_1 -
u_2\|^2_{W_{\tilde g_{R,\phi}} \oplus C^{2,\alpha}_{-\delta}(M)} \\ 
& \leq & 2C_Q \|G_{R,\phi}\| \max\{\|u_1\|, \|u_2\| \} \|u_1 - u_2\| \\ 
& \leq & \frac{1}{2} \|u_1 - u_2\|.
\end {eqnarray*}
\hfill $\blacksquare$

If $u = w + v$ is the unique fixed point of $K_{R,\phi}$ then we
define $g_{R,\phi}$ to the solution to our nonlinear gluing problem:
$$g_{R,\phi} = (1 + v)^\frac{4}{n-2} \tilde g_{R,\phi} (w).$$

\section {Nondegeneracy of the Solution}

In this section we will first prove some preliminary lemmas and then
show that for $R$ sufficiently large the metric $g_{R,\phi}$ is
unmarked nondegenerate. The preliminary lemmas in section \ref{prelim_lemma}
seem to be interesting in their own right. 

\subsection {Some Preliminary Lemmas} \label {prelim_lemma}

In order to prove that $g_{R,\phi}$ is unmarked nondegenerate for sufficiently
large $R$ we will need the following lemmas.

\begin {lemma} \label {balancing}
Suppose $v \in B_{g_1}$ decays like $e^{-\delta_j t_j}$ near all $p_j$
for some $\delta_j > 1$ except
$p_0$. Then $v \sim a v_\epsilon^{0,+}$ near $p_0$ for some $a \in
\mathbb{R}$. Moreover, if there exists $w \in B_{g_1}$ with  $w \sim
v_\epsilon^{0,-}$ near $p_0$ then $v$ decays at least like $e^{-\delta
  t_0}$ for some $\delta > 1$ near $p_0$.  
\end {lemma}
{\bf Proof}: We know that $v \sim \sum (\alpha_{i,+} v_\epsilon^{i,+}
+ \alpha_{i,-} v_\epsilon^{i,-}) \eta_i$
near $p_0$, so suppose $\alpha_{0,-} \neq 0$. Recall from remark
\ref{translations} that we have a Jacobi field $w \in B_{g_1}$ such
that $w \sim v_\epsilon^{0,+}$ near $p_0$ and such that $w \in
C^{2,\alpha}_1(M_1)$. Then 
\begin {eqnarray*} 
0 & = & \lim_{r \rightarrow 0} \int_{M_1 \backslash (\cup B_r(p_j))} w
L_{g_1} v - v L_{g_1}w \\ 
& = & \lim_{r \rightarrow 0} \int_{\partial M_1 \backslash (\cup
B_r(p_j))} w \frac{\partial v}{\partial \nu} - v \frac{\partial
w}{\partial \nu} \\ 
& = & \lim_{r \rightarrow 0} [\int_{t_0 = -\log r} (v_\epsilon^{0,+}
(v_\epsilon^{0,-})' - v_\epsilon^{0,-} (v_\epsilon^{0,+})')  +
O(r^{\gamma_{n+1}(\epsilon)}) + O(\sup_{j\geq 1} r^{\delta_j})]. 
\end {eqnarray*} 
But this last term is just the Wronskian of $v_\epsilon^{0,+}$ and
$v_\epsilon^{0,-}$, which one can write explicitly as $v_\epsilon^{0,+}
(v_\epsilon^{0,-})' - v_\epsilon^{0,-} (v_\epsilon^{0,+})' =  e^{u_\epsilon^2
- u^2_\epsilon (0)}$, which is bounded away from zero. The proof of the
remaining case of this lemma uses an identical argument. \hfill $\blacksquare$

We will also need a lemma regarding Delaunay metrics on finite
cylinders. Before we can state this lemma we need to define the
following function space.
\begin {defn}
The function space $C^{l,\alpha}_\delta([-T,T] \times S^{n-1})$ is
defined to be the space of functions such that 
$$\|u\|_{C^{l,\alpha}_\delta} = \sup_{|\hat t| \leq T-1} \|\left (
\frac{\cosh^\delta(T)}{\cosh^\delta (t)} \right ) \cdot
u\|_{C^{l,\alpha}([\hat t -1, \hat t + 1] \times S^{n-1})}$$ is finite.
\end {defn}

\begin {lemma} \label {finite-cyl}
Pick $\delta \in (1, \gamma_{n+1}(\epsilon))$. Then there exists an operator 
$$H_{T,\epsilon} : C^{0,\alpha}_{-\delta}([-T,T]\times S^{n-1})
\rightarrow C^{2,\alpha}_{-\delta}([-T,T] \times S^{n-1}) $$ 
such that $u = H_{T,\epsilon}(f)$ solves $L_{g_\epsilon}(u) = f$
Moreover, $H_{T,\epsilon}$ is uniformly bounded in $T$.  
\end {lemma} 

Notice we do not say anything about the boundary values of
$H_{T,\epsilon}(f)$, other than that they are bounded by $\|H\| \cdot
\|f\|_{C^{2,\alpha}_{-\delta}([-T,T] \times S^{n-1})}$. 

{\bf Proof}: First choose a cutoff function $\beta$ on $[-T,T] \times
S^{n-1}$ such that 
$$\beta(t,\theta) = \left \{ \begin {array}{rl} 1 & t \leq -1 \\ 0 & t
\geq 1. \end {array} \right .$$
Next let $u_1 \in C^{2,\alpha}_{-\delta} ([-T,\infty) \times S^{n-1})$
solve $L_{g_\epsilon}(u_1) = \beta f$. We can find such a solution
because $g_\epsilon$ is unmarked nondegenerate and we can use
$v_\epsilon^{i,\pm}$ to eliminate the part of $u_1$ which grows at a
rate of $e^t$ or less. Similarly let $u_2 \in C^{2,\alpha}_{\delta}
((-\infty, T]\times S^{n-1})$ solve $L_{g_\epsilon} (u_2) =
(1-\beta)f$. If we let $\hat H_{T,\epsilon}(f) 
= \beta u_1 + (1-\beta) u_2$ then the lemma follows from a
perturbation argument as in the proof of proposition \ref{nondegen}. 

\hfill $\blacksquare$

\subsection {The Nondegeneracy}

In this section we complete the proof of theorem \ref {main-thm} by
showing that $g_{R,\phi}$ is unmarked nondegenerate for $R$ sufficiently
large. We will argue by contradiction, assuming that for some sequence
$R_l \rightarrow \infty$ the metrics $g_l = g_{R_l, \phi}$ are
unmarked degenerate. 

Thus we can find $\delta_l > 1$ and  $0 \neq u_l \in
C^{2,\alpha}_{-\delta_l} (M)$ such that $L_{g_l} (u_l) = 0$. We will
normalize $u_l$ so that 
$$\sup_M \rho_l^{-1} |u_l| = 1$$ where $\rho_l$ is a positive
weighting function we will define in the next paragraph.

First choose for $\delta \in (1, \inf \{ \gamma_{n+1}(\epsilon_j),
\gamma_{n+1}(\epsilon_j')\})$ and recall that for each $l$ we can
decompose $M$ as  
$$M = M_1^c \cup M_2^c \cup (\cup_1^{k_1-1} B_{r_0}(p_j) \backslash
\{p_j\}) \cup (\cup_1^{k_2-1} B_{r_0}(q_j) \backslash \{q_j\}) \cup
\hat C_{R_l}.$$ Then we define the weighting function $\rho_l$ by 
$$\rho_l(p) = \left \{ \begin {array} {rl} 1 & p \in M_i^c \\
e^{-\delta t_j} & p \in B_{r_0/2}(p_j) \backslash \{p_j\} 
\mbox{ for } j = 1 \dots k_1-1\\  
e^{-\delta \tau_j} & p \in B_{r_0/2}(q_j) \backslash \{q_j\} 
\mbox{ for } j = 1 \dots k_2-1\\
\frac{\cosh^\delta R_l}{\cosh^\delta t} & p = (t, \theta) \in \hat C_{R_l}
\end {array} \right.$$
Let $p_l$ be a point where the supremum is achieved. Notice we always
have $|u_l(p)| \leq \rho_l(p)$, with equality at $p_l$. We will obtain
various contradictions depending on where $p_l$ occurs.

First consider the case where $p_l = (t_l, \theta_l) \in \hat C_{R_l}$
with $|t_l|$ bounded. In this case we restrict to $\hat C_{R_l}$ and
renormalize by setting 
$$\tilde u_l (t,\theta) = (\cosh^{-\delta} R_l) u_l(t,\theta).$$
Then choose a subsequence which converges uniformly on compact sets
and such that $(t_l, \theta_l) \rightarrow (\bar t, \bar \theta)$. In
the limit we obtain a Jacobi field $\bar u$ for the Delaunay metric
$g_\epsilon$ such that $$|\bar u(t, \theta)| \leq \cosh^{-\delta} t$$
with equality at $(\bar t, \bar \theta)$, which is a contradiction.

Next consider the case where $p_l = (t_l, \theta_l) \in C_{R_l}$ with
$|t_l|$ and $|t_l \pm R_l|$ all unbounded. We will treat the instance
where $t_l <0$; the case where $t_l >0$ is similar. In this case we
restrict to the part of $\hat C_{R_L}$ parameterized by $(t,\theta) \in
[-R_l -t_l, |t_l|] \times S^{n-1}$ and renormalize by setting 
$$\tilde u_l(t,\theta) = \left ( \frac{\cosh^\delta t_l}{\cosh^\delta
R_l} \right ) u_l(t + t_l, \theta).$$ With this renormalization
$|\tilde u_l(0, \theta_l)| = 1$ and 
\begin {eqnarray*}
|u_l(t,\theta)| & \leq & \frac{\cosh ^\delta t_l}{\cosh^\delta
 (t+t_l)} \\
& \leq & 2^\delta (e^t + e^{-t -2t_l})^{-\delta} \leq 2^\delta e^{\delta
 t}.
\end {eqnarray*}
Extracting a convergent subsequence we obtain a Jacobi field $\bar u$
for $g_\epsilon$ such that $|\bar u(0, \bar \theta)| = 1$ and $|\bar
u(t, \theta)| \leq 2^\delta e^{\delta t}$, which is a contradiction. 

Next consider the case where $p_l = (t_{j,l}, \theta_{j,l}) \in
B_{r_0}(p_j) \backslash \{p_j\}$ (for $j = 1 \dots k_1-1$) and $t_l
\rightarrow \infty$. In this case we restrict to $B_{r_0}(p_j)
\backslash \{p_j\}$, renormalize by setting 
$$\tilde u_l(t_j, \theta_j) = e^{\delta t_{j,l}} u_l(t_j + t_{j,l},
\theta_j)$$ and argue as in the previous case. The case where $p_l =
(\tau_{j,l}, \theta_l) \in B_{r_0}(q_j) \backslash \{q_j\}$ with
$\tau_{j,l} \rightarrow \infty$ is similar. 

Next consider the case where $p_l \in \Omega_1$, where $\Omega_1$ is
some fixed compact set containing $M_1^c$. Notice $\rho_l$ is bounded
and bounded away from $0$ in $\Omega_1$. Restrict to 
$$M_1^c \cup (\cup_1^{k_1-1} B_{r_0}(p_j) \backslash \{p_j\}) \cup
\{(t, \theta) \in \hat C_{R_l} : t < 0\}$$ and take a subsequence which
converges uniformly on compact sets (and so $p_l \rightarrow \bar p
\in \Omega_1$). Then in the limit we obtain a Jacobi field $\bar u$ on
$M_1$ which decays exponentially near all $p_j$ except $p_0$, and
which has subexponential growth near $p_0$. Also, by the normalization
$|u(\bar p)| \neq 0$. By lemma \ref{balancing}, $\bar u$  must also
decay like $e^{-\delta t_0}$ near $p_0$ for some $\delta > 1$, which
contradicts the unmarked nondegeneracy of $(M_1, g_1)$. 

Finally, consider the case where $p_l \in \Omega_2$, where $\Omega_2$
is a fixed compact set containing $M_2^c$. If we restrict to $$M_2^c
\cup (\cup_1^{k_2-1} B_{r_0}(q_j) \backslash \{q_j\}) \cup
\{(t,\theta) \in \hat C_{R_l} : t>0\}$$ and choose a convergent subsequence
as we did in the previous case, we can only conclude that the limit
$\bar u$ is asymptotic to $a v_\epsilon^{0,+}$ near $q_0$. At this
point we rescale so that $a=1$. Fix some $R_0>0$. Then there is an
$l_0$ depending on $R_0$ such that for $l \geq l_0$ 
$$\|u_l - v_\epsilon^{0,+}\|_{C^{2,\alpha}([R_l-R_0-1, R_l] \times S^{n-1})}
= O(e^{-\delta R_0}).$$
Moreover, if we restrict $u_l$ to 
$$M_1^c \cup (\cup_1^{k_1-1} B_{r_0}(p_j) \backslash \{p_j\}) \cup
\{(t, \theta) \in \hat C_{R_l} : t < 0\}$$ we know by the previous argument
that it must converge uniformly to zero. Thus for $l \geq l_0$
$$\|u_l\|_{C^{2,\alpha}([-R_l,-R_l+R_0 +1] \times S^{n-1})} = O(e^{-\delta
R_0}).$$ 

Recall that we can write the metric $g_l = g_{R_l,\phi}$ on $\hat
C_{R_l}$ as  
$$(u_\epsilon + v_l)^\frac{4}{n-2}(dt^2 + d\theta^2)$$ where 
$$|v_l(t,\theta)| = O\left ( \frac{\cosh^{\gamma_{n+1}(\epsilon)}
  t}{\cosh^{\gamma_{n+1}(\epsilon)} R_l} \right ) .$$ 
Thus $L_{g_l} - L_{g_\epsilon}$, when restricted to $\hat
C_{R_l}$ a second order differential operator whose coefficients are
$O(\frac{\cosh^\mu t}{\cosh^\mu R_l})$ on $\hat C_{R_l}$ for some $\mu
\in (\delta, \gamma_{n+1}(\epsilon))$, which implies 
$$L_{g_\epsilon}(u_l)(t,\theta) = O \left ( \frac{\cosh^{\mu - \delta} 
t}{\cosh^{\mu - \delta} R_l} \right )$$ on $\hat C_{R_l}$. In other words,
$L_{g_\epsilon}(u_l) \in C^{0,\alpha}_{\delta - \mu}([-R_l+R_0,R_l-R_0]
\times S^{n-1})$ and 
$$\|L_{g_\epsilon}(u_l)\|_{C^{0,\alpha}_{\delta - \mu}([-R_l+R_0, R_l-R_0]
\times S^{n-1})} = O(e^{(\delta - \mu) R_0}).$$ Let  
$$\tilde u_l = H_{R_l-R_0, \epsilon}(L_{g_\epsilon} (u_l)).$$
Then 
\begin{eqnarray*} 
0 & = & \int_{[-R_l+R_0, R_l-R_0] \times S^{n-1}} \tilde u_l
L_{g_\epsilon}(v_\epsilon^{0,-}) - v_\epsilon^{0,-} L_{g_\epsilon}(\tilde
u_l) \\ 
& = & \int_{\{R_l-R_0\}\times S^{n-1}} (\tilde u_l \frac{\partial
v_\epsilon^{0,-}}{\partial \nu} - v_\epsilon^{0,-} \frac{\partial \tilde
u_l}{\partial \nu}) - \int_{\{-R_l+R_0\}\times S^{n-1}} (\tilde u_l
\frac{\partial v_\epsilon^{0,-}}{\partial \nu} - v_\epsilon^{0,-}
\frac{\partial \tilde u_l}{\partial \nu}) \\ 
& = & \int_{\{R_l-R_0\}\times S^{n-1}} (v_\epsilon^{0,+} \frac{\partial
v_\epsilon^{0,-}}{\partial t} - v_\epsilon^{0,-} \frac{\partial
v_\epsilon^{0,+}}{\partial t} + O(e^{-\delta R_0} + e^{(\delta - \mu)
R_0}))  + \int _{\{-R_l+R_0\} \times S^{n-1}} O(e^{-\delta R_0})\\ 
& = & 1 + O(e^{-\delta R_0} + e^{(\delta - \mu) R_0}).  
\end {eqnarray*}
This completes the proof of theorem \ref{main-thm}.

\hfill $\blacksquare$

\section {Questions}

In this final section we raise some interesting questions related to
this construction.

The first question is: how much of this gluing construction can be
extended to arbitrary (complete, connected and noncompact) manifolds
with constant positive scalar curvature? The first requirement we see
is that all the ends must be asymptotically Delaunay. So it might be
natural to apply this theorem with $M_i = \bar M_i \backslash \Lambda$
where $\bar M_i$ is a closed locally conformally flat manifold and
$\Lambda$ is a finite set. In this case the ends of $M_i$ correspond
to punctured neighborhoods of $p \in \Lambda$ and the metric is indeed 
asymptotically Delaunay there. However, we lack the conformal Killing
fields which give rise to the asymptotic translations in this case (as
we used in remark \ref {translations}). However, one might be able to
prove a similar result where one supposes that the localization of
$B_{g_1}$ to $E_o$ is trivial, i.e. that there are no bounded Jacobi
fields for $g_1$ which decay on all ends but $E_0$. This may not be
such a strong hypothesis, because as argument similar to that of lemma
\ref{balancing} shows that the localization of $B_{g_1}$ to $E_0$ can
be at most $1$-dimensional. 

The second question stems from conversations with N. Korevaar and is:
how much can we say about the global structure 
of the moduli space of complete scalar curvature metrics $n(n-1)$
metrics on $S^n \backslash \{p_1 \dots p_k\}$? The simplest nontrivial
case seems to be $k=3$. We can conclude (using an Alexandrov
reflection argument as in \cite {CGS} and \cite {KKS}) that these
metrics must be symmetric under reflection through some equatorial
$S^{n-1}$ (after composing with a conformal motion of $S^n$). More
precisely, one can use stereographic projection to turn the problem
into a scalar PDE on $\mathbb{R}^n$ and take inversion through
$n-1$-dimensional spheres centered at the origin in place of
reflection through some hyperplane. The same arguments as in \cite
{KKS} hold. However, the main tool we lack in the scalar curvature
case is a way to get necksize bounds on the ends, and reconstruct the
metric from the asymptotic necksizes as in \cite
{GKS}. Grosse-Brauckmann, Kusner and Sullivan use a conjugate minimal
surface in $S^2$ to classify all three ended, genus zero constant mean
curvature surfaces in \cite {GKS}. At the present, we do not have any
way to either find necksize bounds or show that all possible
combination of necksizes allowed by balancing is realized.

\end {document}